\documentclass[12pt]{article}

\newcommand{\eop}{\hfill $/$\hspace*{-.1cm}$/$\hspace*{-.1cm}$/$\vspace{.1in}}

\newtheorem{lemma}{Lemma}

\newtheorem{corollary}[lemma]{Corollary}

\newtheorem{proposition}[lemma]{Proposition}

\newcommand{\tworightarrows}{\stackrel{\displaystyle \rightarrow}{\rightarrow}}

\begin{document}

\section*{Calculating maps between $n$-categories}

Carlos Simpson
\newline
CNRS, Universit\'e de Nice-Sophia Antipolis, 
Parc Valrose
\newline
06108 Nice Cedex 2,
France

\begin{center}
{\em Introduction}
\end{center}

This short note addresses the problem of how to calculate in a
reasonable way the homotopy classes of maps between two $n$-categories
(by which we mean Tamsamani's $n$-nerves \cite{Tamsamani}). The closed model structure of
\cite{cmstruct} gives an abstract way of calculating this but it isn't
very concrete so we would like a more down-to-earth calculation.  For
the purposes of the present note we shall use the closed model structure
of \cite{cmstruct} to prove that our method gives the right answer. It
might be possible to develop the theory entirely  using the present
calculation as the definition of ``map'', but that is left as an open 
problem.

Denote by $n-Cat$ the category of $n$-categories introduced by
Tamsamani
\cite{Tamsamani}, and by $L(n-Cat)$ the Dwyer-Kan simplicial
localization \cite{dk} where we divide out by the equivalences between
$n$-categories (Tamsamani calls these ``external equivalences'').
Note that the $1$-category obtained from 
$L(n-Cat)$ by applying $\pi _0$ to the simplicial $Hom$ sets
is just the Gabriel-Zisman \cite{gz} localization $Ho(n-Cat)$ which 
was introduced
in \cite{Tamsamani}.

We work with the model category $nPC$ of $n$-precats introduced in
\cite{cmstruct}, and we shall adopt the notations as from there.
In particular, $nPC$ is the category of presheaves over the category
$\Theta ^n$.

In \cite{cmstruct} the procedure for calculating the homotopy classes
of maps from $A$ to $B$ consists of choosing a fibrant replacement $B\rightarrow
B'$ and looking at maps from $A$ to $B'$; two maps are homotopic if there is a
homotopy $A\times \overline{I}\rightarrow B$ relating them (where $\overline{I}$
is the $1$-category with two isomorphic objects). The problem is that the
notion of fibrant replacement is not very explicit, depending on the addition of
a wide class of pushouts by trivial cofibrations to $B$ to get $B'$. Here we 
would like to pursue the alternative strategy of replacing $A$ by a 
``free cofibrant''
object $F\rightarrow A$, much as was done for morphisms
of diagrams in Bousfield-Kan \cite{bk}. Then we look at maps from $F$ to $B$.
The advantage is that $F$ can be made
much more explicit, and we shall only need to assume that $B$ is an
$n$-category.

The interest of having such a description developed mainly
out of conversations with A. Hirschowitz during and in the aftermath of
\cite{descente}.

The notion of free cofibration is a basic consequence 
of the shape of cells coming from 
Tamsamani's original ``constancy conditions'' and the corresponding
quotient construction for $\Theta ^n$ as a quotient of $\Delta ^n$.

At the first level in the multisimplicial approach to $n$-categories,
one meets cells which look like $k$-simplices for any $k$. A $k$-simplex
corresponds to $k$ one-morphisms which can be composed, together with their
compositions. On this
level, the notion of free cofibration is easy to imagine: it is just the
inclusion of the set of $k+1$ vertices (i.e. the objects), 
into the $k$-simplex. Our
fundamental observation is that this definition extends into the other
simplicial directions (corresponding to higher arrows) in a natural 
way. 

Once one has noticed the notion of ``free cofibration'' of $n$-precats, the
argument is a straightforward application of postmodern closed model
category techniques as developped principally by Reedy, Dwyer, Kan, 
Hirschhorn \cite{Reedy} \cite{dk} \cite{dhk} \cite{Hirschhorn}.

We end up in Corollary \ref{explicit} with a very explicit description 
of how to calculate the set of homotopy classes of maps from $A$ to $B$ in
$Ho(n-Cat)$.

I would like to thank  C. Berger, A. Hirschowitz, G. Maltsiniotis,
R. Pelissier,
Z. Tamsamani and B. Toen for
discussions surrounding this topic.

\begin{center}
{\em Free cofibrations}
\end{center}

We define a class of morphisms in $nPC$ called {\em free cofibrations}.
These are analogues of the Bousfield-Kan cofibrations in the theory
of diagrams \cite{bk}. 
This class is generated by the following {\em elementary free
cofibrations}. Suppose $M=(m_1,\ldots , m_k)$ is an object in $\Theta
^n$. Then $h(M)\in nPC$, and we define the {\em elementary free
cofibration with target $h(M)$} denoted
$$
\partial (M) \rightarrow h(M)
$$
in the following way. For any other index object $N\in \Theta ^n$, the
elements of $\partial (M)_N$ are defined as those elements $a\in h(M)_N$
corresponding to maps $a: N\rightarrow M$ such that $a$ factors through
some map of the form 
$$ 
(m_1, \ldots , m_{k-1}, 0) \rightarrow (m_1,\ldots , m_k)=M.
$$
In other words, $\partial (M)$ is the union of the images of the $m_k+1$ morphisms
$\hat{M}\rightarrow M$ where $\hat{M}:=(m_1, \ldots , m_{k-1}, 0)$.

An {\em elementary free cofibration} is any cofibration corresponding to
some $M$ as above. A {\em free cofibration} is any cofibration that is
obtained as a sequential (possibly transfinitely sequential) pushout by 
elementary free cofibrations. In other words the class of free
cofibrations is the class generated by the elementary free ones in the
sense of \cite{dhk}.

\begin{lemma}
\label{description}
For any $M$ the ``boundary'' $\partial (M)$ is itself obtained from
$\emptyset$
by a
sequence of pushouts along elementary free cofibrations. 

At the top level, $\partial
(M)$ is obtained by glueing together $m_k+1$ copies of $h(\hat{M})$
along $\partial (\hat{M})$.
\end{lemma}
{\em Proof:}
The statement about what happens at the top level follows directly from
the shapes of Tamsamani's cells; or alternatively from the definition of the
category $\Theta ^n$ as a quotient of $\Delta ^n$. Indeed from the definition,
$\partial (M)$ is clearly the image of the disjoint union $U$ of $m_k+1$ copies of 
$h(\hat{M})$. Two elements of $U_N$ go to the same element of $\partial (M)_N$
if and only if they go to the same element of $h(M)_N$. Thus two elements
of $U_N$ represented by $a,a': N\rightarrow \hat{M}$ (corresponding to different
copies of $h(\hat{M})$ i.e. to different maps $0\rightarrow m_k$)
go to the same element,
if and only if their compositions $N\rightarrow M$ are the same.
In view of the definition of $\Theta ^n$, this holds if and only if $a$ and 
$a'$ factor through some map $N\rightarrow \hat{\hat{M}}$. This gives the
last statement of the lemma.

The first statement follows by induction
in the following way. Let $M^{[i]}$ denote the object $(m_1,\ldots , m_i, 0)$.
Thus if $M$ has length $k$ then $\hat{M}=M^{[k-1]}$. By induction using the
statement about what happens on the top level, we get that 
$\partial (M^{[i]})$ is obtained from $\partial (M^{[i-1]})$ by adding on
$m_i+1$ copies of $h(M^{[i-1]})$. Each of these copies is obtained by taking the
coproduct with the
cofibration
$$
\partial (M^{[i-1]}) \rightarrow h(M^{[i-1]}).
$$
\eop

The main point behind the notion of free cofibrations is the following 
homotopy lifting property.

\begin{proposition}
\label{main}
Suppose $B$ is an $n$-category and $i:B\rightarrow B'$ is a fibrant
replacement.
Then any free cofibration satisfies the up-to-homotopy lifting property
for $i$. More precisely, if $A\rightarrow C$ is a free cofibration and
if 
$$
\begin{array}{ccc}
A & \rightarrow & B \\
\downarrow & & \downarrow \\
C & \rightarrow & B'
\end{array}
$$
is a commutative diagram, then there is a northeast diagonal map
$g: C\rightarrow B$ which
is a lifting up to homotopy, i.e.
there exists a homotopy $C\times
\overline{I}\rightarrow B'$ fixing $A$, starting at the original map
$C\rightarrow B'$ and ending at a map which factors through 
$g:C\rightarrow B$. 
\end{proposition}
{\em Proof:}
It suffices to consider the case where the cofibration $A\rightarrow C$
is of the form $\partial (M)\rightarrow h(M)$. Let $k$ be the length of
$M$
and denote by $m:= m_k$.
A map $h(M)\rightarrow
B$ (resp. $h(M)\rightarrow B'$)
is the same thing as an object of the set
$$
B_{M} = B_{\hat{M},m} \;\;\; (resp.   B'_{M} = B'_{\hat{M},m})
$$
Let $E$ (resp. $E'$) denote the $n+1-k$-category $B_{\hat{M}/}$
(resp. $B_{\hat{M}/}$). Then a map $f:h(M)\rightarrow B$ (resp. $f':
h(M)\rightarrow B'$) is
the same thing as an object of $E_m$ (resp. $E'_m$). 

Fix a morphism $a: \partial (M) \rightarrow B$. This may be thought of
as a collection of elements
$$
\alpha _0,\ldots , \alpha _m \in B_{\hat{M}}
$$
which agree over $\partial (\hat{M})$. In turn these correspond to
objects of $E$. A map $f$ or $f'$ as above which restricts to $a$ on
$\partial (M)$ is just an element 
$$
f\in E_{m/}(\alpha _0,\ldots , \alpha _m)
$$
or
$$
f'\in E'_{m/}(\alpha _0,\ldots , \alpha _m).
$$
The fact that $B\rightarrow B'$ is an equivalence of $n$-categories
yields the statement that the map 
$$
i(a):E_{m/}(\alpha _0,\ldots , \alpha _m)\rightarrow
E'_{m/}(\alpha _0,\ldots , \alpha _m)
$$
is an equivalence of $n-k$-categories \cite{Tamsamani}. 
In particular, given $f'$ in the
right side, there is an $f$ on the left which maps to an object
equivalent to $f'$. It is straightforward to turn this equivalence
between $i(a)f$ and $f'$ into
a homotopy from $f'$ to $i\circ f$, using the fibrant property of $B'$ and
traditional closed model category techniques. 
\eop

The following lemma establishes a relationship between lifting
for the elementary free cofibrations, and being an equivalence.

\begin{lemma}
\label{lifting}
Suppose that $f:A\rightarrow B$ is a morphism of $n$-precats, such that
$B$ is an $n$-category and such that $f$ satisfies the 
on-the-nose lifting property
for every elementary free cofibration. Then $A$ is an $n$-category and
$f$ is an equivalence.
\end{lemma}
{\em Proof:}
We prove this by induction on $n$. If $m\in \Delta$ and
 $N\in \Theta ^{n-1}$ then
a diagram 
$$
\begin{array}{ccc}
\partial (m,N) & \rightarrow & h(m,N)\\
\downarrow && \downarrow \\
A & \rightarrow & B
\end{array} .
$$ 
is the same thing as a choice of $m+1$ objects $x_0,\ldots , x_m\in A_0$
together with a diagram
$$
\begin{array}{ccc}
\partial (N) & \rightarrow & h(N)\\
\downarrow && \downarrow \\
A_{m/}(x_0,\ldots , x_m) & \rightarrow & B_{m/}(f(x_0),\ldots , f(x_m)).
\end{array}
$$ 
This can be seen by using the description given by Lemma \ref{description}.

Similarly, a lifting $h(m,N)\rightarrow A$ is the same thing as a
lifting $h(N)\rightarrow A_{m/}(x_0,\ldots , x_m)$.

In particular if $f$ satisfies the lifting property then so do the
morphisms 
$$
f_{m/}(x):A_{m/} (x_0,\ldots , x_m) \rightarrow  B_{m/} (fx_0,\ldots ,
fx_m).
$$
By induction we conclude that the 
$A_{m/} (x_0,\ldots , x_m)$ are $n-1$-categories and that these
morphisms $f_{m/}(x)$ are equivalences. In turn, the fact that $B$ is an
$n$-category and that these morphisms are equivalences, implies that $A$
is an $n$-category and that $f$ is fully faithful. Essential
surjectivity of $f$ follows from surjectivity on objects which is a limiting case
of the lifting property.
\eop

{\em Remark:} It would probably be a good idea to replace the notion of 
``easy equivalence'' which was used at the start of \cite{cmstruct}, by
the notion of morphism satisfying the lifting property of this lemma.

\begin{center}
{\em Good resolutions}
\end{center}

Next we apply the strategy of Reedy-Dwyer-Kan for calculating the space
of maps from an $n$-category $A$ to another one $B$ in the localized
simplicial category $L(n-Cat)$, by using a cosimplicial resolution of
$A$. The trick is to use a cosimplicial resolution which is Reedy-cofibrant in
terms of free cofibrations in $nPC$.

Recall that a cosimplicial object $F^{\cdot}$ in a model category $M$ is
{\em Reedy-cofibrant} if for every $k$ the morphism
$$
Latch ^k(F) \rightarrow F^k
$$
is a cofibration in $M$. In our situation, we are not sure whether there
exists a closed model structure on $nPC$ which has free cofibrations as
its fibrations (there probably does...). However, it still makes sense
to ask that a cosimplicial object $F^{\cdot}$ be Reedy-cofibrant with
respect to the class of free cofibrations: this just means that we ask
for every $k$ for the morphism 
$$
Latch^k(F) \rightarrow F^k
$$
to be a free cofibration. We call such a resolution, {\em Reedy-free-cofibrant}.

The cornerstone of our calculation is the construction of a
Reedy-free-cofibrant resolution of any $n$-category $A$.

\begin{proposition}
\label{cornerstone}
Suppose $A$ is an $n$-category. Then there is a natural functorial
cosimplicial object $F^{\cdot}$ in $nPC$ together with a morphism to the
constant cosimplicial object associated to $A$ (we denote this by 
$F^{\cdot} \rightarrow A$) such that each stage $F^k\rightarrow A$ is an
equivalence
of $n$-categories, and such that $F^{\cdot}$ is Reedy-cofibrant with
respect to the class of free cofibrations in $nPC$.
\end{proposition}
{\em Proof:}
There are two canonical ways to complete any morphism of
$n$-precats
$U\rightarrow V$ to a diagram 
$$
U \rightarrow W \rightarrow V
$$
where the first  morphism is a free cofibration and the second morphism
satisfies the lifting property with respect to the elementary free
cofibrations. 

The ``small'' way is to complete the diagram by adding in lifts of all
elementary free cofibrations which don't already have lifts. This is of course
to be prefered in calculations, but it is not functorial. 

The ``big'' way to complete the above diagram is to add in
lifts of all elementary free cofibrations (without looking to see whether there
is already a lift), and iterate over the first infinite
ordinal $\omega$. The big way is functorial, so we adopt it for the theory.

If $V$ is an $n$-category then by Lemma \ref{lifting}, $W$ is also an
$n$-category and the morphism $W\rightarrow V$ is an equivalence.

Using this construction, the standard construction of the
Reedy-cofibrant cosimplicial resolution (adding in each stage by
enforcing the lifting property) works to give the desired
$F^{\cdot}$. 

In order to deal with the degeneracies (although it isn't clear to me that we
really need them) one needs the following observation: if $U\rightarrow
V\leftarrow W$ is a diagram of $n$-categories such that 
both maps satisfy the lifting property of Lemma \ref{lifting}, then the fiber
product $U\times _VW\rightarrow V$ also satisfies the lifting property and 
in particular it is again an $n$-category equivalent to the other ones. 

We explicitly describe the first couple of steps in the standard construction.
First choose a free cofibrant $F^0$ with map $F^0\rightarrow A$ satisfying
the lifting property of Lemma \ref{lifting}. Then choose $F^1$ fitting into 
the diagram 
$$
F^0\sqcup F^0 \rightarrow F^1\rightarrow F^0
$$ 
so that the first map is a free cofibration and the second map 
(the ``degeneracy'') satisfies
the lifting property. Let $C$ be the pushout
$$
C:= (F^1\cup ^{F^0}F^1) \cup ^{F^0\sqcup F^0} F^1,
$$
in other words it is three copies of $F^1$ glued together at three copies of
$F^0$ in the form of a triangle. It is the latching object: $C=Latch^2(F)$.
We have a map
$$
C\rightarrow F^1\times _{F^0}F^1
$$
where the first projection is the identity in the first and last components of
$C$
and the degeneracy on the middle component; and the second projection is the
identity on the second and third components and the degeneracy on the first
component
(for the author at least it was easier to write down the dual map for
simplicial objects).
From the previous paragraph the target 
$F^1\times _{F^0}F^1$ is again an $n$-category mapping
via an equivalence to $F^0$. Choose $F^2$ to fit into a diagram
$$
C\rightarrow F^2\rightarrow  F^1\times _{F^0}F^1
$$
where the first map is a free cofibration and the second map satisfies the
lifting property of Lemma \ref{lifting}. 

For the rest of the cosimplicial object, continue in the same way.

The whole cosimplicial object projects to $F^0$ by the degeneracy maps, and
this in turn maps to $A$. We obtain the map $F^{\cdot}\rightarrow A$
which on every stage is an equivalence of $n$-categories $F^k\cong A$. 
\eop

Applying the up-to-homotopy lifting property
of Proposition \ref{main} to the resolutions provided
by Proposition \ref{cornerstone} will give the calculation we are
looking for. This is resumed in the following proposition.

\begin{proposition}
\label{calculation}
Suppose $F^{\cdot}\rightarrow A$ is a Reedy-free-cofibrant cosimplicial
resolution
as in Proposition \ref{cornerstone}. Suppose $i:B\rightarrow B'$ is a
fibrant replacement of an $n$-category $B$. Then the morphism of
simplicial sets
$$
Hom _{nPC}(F^{\cdot}, B) \rightarrow
Hom _{nPC}(F^{\cdot}, B')
$$
is a weak equivalence. 
\end{proposition}
{\em Proof:}
This general property of maps having the up-to-homotopy lifting property
with respect to a class of ``cofibrations'' probably goes back to Reedy
\cite{Reedy}. We give an argument for completeness.

Suppose $K$ is a finite simplicial set. 
Recall that by adjunction we can define an operation 
$$
K\mapsto K\otimes F^{\cdot} \in nPC
$$
such that a
map
$$
K\rightarrow Hom_{nPC} (F^{\cdot}, B)
$$
is the same thing as a map $K\otimes F^{\cdot}\rightarrow B$.
We claim that if $K\rightarrow L$ is a cofibration of finite
simplicial sets then
$$
K\otimes F^{\cdot} \rightarrow L\otimes F^{\cdot}
$$
is a free cofibration in $nPC$.  Indeed, for the generating cofibrations
of simplicial sets $\partial \Delta ^k \rightarrow \Delta ^k$ we just
get back the original free cofibrations $Latch^k(F)\rightarrow F^k$ in
the definition of Reedy-free-cofibrant resolution. 

Suppose $(\ast )$
$$
\begin{array}{ccc}
K & \rightarrow & Hom_{nPC} (F^{\cdot}, B) \\
\downarrow && \downarrow \\
L & \rightarrow &  Hom_{nPC} (F^{\cdot}, B')
\end{array} 
$$  
is a diagram. It corresponds to a diagram
$$
\begin{array}{ccc}
K\otimes F^{\cdot} & \rightarrow & B \\
\downarrow && \downarrow \\
L\otimes F^{\cdot} & \rightarrow &  B'
\end{array} 
$$  
and since we saw above that the left vertical arrow is a free
cofibration, Proposition \ref{main} gives a homotopy lifting. This means that
there is a lifting in the original diagram such that the two resulting
morphisms
$$
L \tworightarrows Hom_{nPC} (F^{\cdot}, B')
$$
come from maps 
$$
L\otimes F^{\cdot}  \tworightarrows   B'
$$
which are related by a homotopy
$$
L\otimes F^{\cdot} \times \overline{I}  \rightarrow   B'
$$
or equivalently 
$$
L\otimes F^{\cdot}\rightarrow \underline{Hom} (\overline{I} ,   B').
$$
Thus the two maps $L \tworightarrows Hom_{nPC} (F^{\cdot}, B')$
are related by a ``homotopy''  which is a map
$$
L \rightarrow Hom_{nPC} (F^{\cdot}, \underline{Hom} (\overline{I} ,  
B')).
$$
The operation $ Hom_{nPC} (F^{\cdot}, -)$ is homotopy-invariant (and
takes fibrations to fibrations) when the
argument is a fibrant object of $nPC$. Thus it takes a path-space object
in $nPC$ in the sense of Quillen, such as $\underline{Hom}(\overline{I}, B')$,
to a path-space object in the Kan
closed model category of simplicial sets. In particular, if two maps
$$
L \tworightarrows Hom_{nPC} (F^{\cdot}, B')
$$
are related by a ``homotopy'' in the above sense, then they are related
by a Quillen homotopy in the closed model category of simplicial sets.
We have now shown that the morphism of simplicial sets 
$$
Hom _{nPC}(F^{\cdot}, B) \rightarrow
Hom _{nPC}(F^{\cdot}, B')
$$
satisfies the up-to-homotopy lifting property with respect to any
cofibration 
$K\rightarrow L$ of finite simplicial sets. This implies that it is a
weak equivalence.
\eop

Note in the situation of the above proposition, the $F^k$ and $B$ are
both $n$-categories. Recall that $n-Cat$ is a full subcategory of $nPC$
so the morphism which is an equivalence by the proposition may be
written as
$$
Hom _{n-Cat}(F^{\cdot}, B) \rightarrow
Hom _{n-Cat}(F^{\cdot}, B').
$$

\begin{corollary}
\label{calculation2}
Suppose $F^{\cdot}\rightarrow A$ is a free-cofibrant cosimplicial
resolution
as in Proposition \ref{cornerstone}, and suppose $B$ is an $n$-category.
Then the simplicial set $Hom _{n-Cat}(F^{\cdot}, B)$ is naturally
equivalent to the simplicial set of morphisms from $A$ to $B$ in the simplicial category
$L(n-Cat)$.
\end{corollary}
{\em Proof:}
In \cite{dk} (see also \cite{dhk}) it is shown that the simplicial set
of morphisms in the Dwyer-Kan localization is naturally equivalent to
the simplicial set obtained by using a Reedy-cofibrant cosimplicial
resolution of the domain object and assuming that the range-object is
fibrant. The resolution $F^{\cdot}\rightarrow A$ is Reedy-cofibrant for
the usual notion of cofibration in the closed model structure of
\cite{cmstruct},
because free cofibrations are in particular cofibrations. Thus 
there is a natural equivalence of simplicial sets
$$
Hom _{L(nPC), \cdot}(A,B) \cong Hom _{nPC}(F^{\cdot}, B').
$$
On the other hand, $L(n-Cat)$ is sandwiched in between $L(nPC_f)$ and
$L(nPC)$ but the latter two are equivalent \cite{dk} so 
$$
Hom _{L(n-Cat), \cdot}(A,B) \cong Hom _{L(nPC), \cdot}(A,B) .
$$
Finally, by the previous proposition and composing, we obtain the equivalence
$$
Hom _{L(n-Cat), \cdot}(A,B) \cong Hom _{nPC}(F^{\cdot}, B)=Hom _{n-Cat}(F^{\cdot}, B).
$$
\eop

\begin{corollary}
\label{calculation3}
Suppose $F^{\cdot}\rightarrow A$ is a free-cofibrant cosimplicial
resolution
as in Proposition \ref{cornerstone}, and suppose $B$ is an $n$-category.
Then the set of morphisms from $A$ to $B$ in $Ho(n-Cat)$ is calculated
as
$$
Hom _{Ho(n-Cat)}(A,B) = \pi _0 Hom _{n-Cat}(F^{\cdot}, B).
$$
\end{corollary}
{\em Proof:}
This follows immediately from the fact that $Ho(n-Cat)$ is just $\pi _0$
of $L(n-Cat)$ (cf \cite{dk} for this fundamental property of the
Dwyer-Kan
localisation $L(-)$).
\eop

\pagebreak[5]

\begin{center}
{\em Further remarks}
\end{center}

In order to calculate the $\pi _0$ (i.e. the set of morphisms in
$Ho(n-Cat)$) we only need the first stage of the Reedy-free-cofibrant
resolution:
$$
F^0 \tworightarrows F^1 \rightarrow A.
$$
This yields a diagram of sets
$$
Hom (F^1, B) \tworightarrows Hom (F^0,B)
$$
and $Hom _{Ho(n-Cat)}(A,B)$ is the set-theoretic coequalizer of these
two morphisms. In particular, from this description two maps
$F^0\rightarrow B$ could represent homotopic maps from $A$ to $B$, but
only be related by a chain of two or more maps $F^1\rightarrow B$.
We shall now show that actually there is no need to refer to chains of
equivalences.

\begin{lemma}
\label{easy}
Suppose $F^0$ is free-cofibrant, and 
$$
F^0 \sqcup F^0 \rightarrow F^1
$$
is a free cofibration, fitting into a diagram
$$
F^0 \tworightarrows F^1 \rightarrow A
$$
where all maps to $A$ are equivalences of $n$-categories.
Then for any $n$-category $B$,
two maps $F^0\rightarrow B$ are homotopic if and only if
their disjoint union extends to a map $F^1\rightarrow B$. 
\end{lemma}
{\em Proof:}
Fix the diagram $F^0 \tworightarrows F^1 \rightarrow A$ in question.
We call a {\em basic homotopy} between maps $F^0\rightarrow B$,
any one which is obtained by a morphism $F^1\rightarrow B$. To prove the
lemma, it suffices to show that the composition of two basic homotopies
is again homotopic to a basic homotopy. For this choose an extension of
the resolution to level $2$, i.e. with $F^2\rightarrow A$ restricting to
three copies of the map $F^1\rightarrow A$. Let $D$ denote the union of
the first two copies of $F^1$ (it is in $nPC$ but not an $n$-category). 
Then $D\rightarrow F^2$ is a weak equivalence of $n$-precats. It is also
a free cofibration. Two
homotopies result in a map $D\rightarrow B$; compose this to obtain a
map $D\rightarrow B'$ which extends to $h':F^2\rightarrow B'$ since $B'$ is
fibrant.
The fact that $D\rightarrow F^2$ is a free cofibration allows us to
apply Proposition \ref{main} to move $h'$ back to a map $h: F^2\rightarrow
B$,
which gives (by restricting to the third copy of $F^1$) a basic homotopy
composing the first two given ones. 
\eop

\begin{corollary}
\label{explicit}
Suppose $A$ is an $n$-category. Choose a diagram 
$$
F^0 \tworightarrows F^1 \rightarrow A
$$
with 
$F^0$ a free cofibrant object and $F^0 \sqcup F^0 \rightarrow F^1$
a free cofibration, and all objects being $n$-categories mapping by
equivalences to $A$. Then for any $n$-category $B$ the set of homotopy
classes of maps from $A$ to $B$ in $Ho(n-Cat)$ may be described as the
set of maps $F^0\rightarrow B$ modulo the relation that two such maps
are equivalent if and only if there exists a map $F^1\rightarrow B$ 
restricting to
their disjoint union on $F^0 \sqcup F^0$.
\end{corollary}
{\em Proof:}
This is what we just showed in the previous lemma.
\eop

{\bf Exercise 1:} {\em Let $C$ denote the $1$-category with three objects
$0,1,2$ and morphisms $f: 0\rightarrow 1$ and $g: 1\rightarrow 2$
composing to $gf: 0\rightarrow 2$. Consider $C$ as a $2$-category, and
calculate a free-cofibrant $2$-category $F$ with an equivalence
$F\rightarrow C$. }

{\bf Exercise 2:} {\em If $F$ is a free cofibrant object then the set of maps
from $F$ to an $n$-category $B$ has an easy inductive description in
terms of the cells which were used to construct $F$. Write this out in detail.}


\begin{thebibliography}{MM}

\bibitem{bk}
A. Bousfield, D. Kan. {\em Homotopy limits, completions and localisations}.
Lecture Notes in Math. {\bf 304}, Springer-Verlag (1972).


\bibitem{dk} W. Dwyer, D. Kan.
\newline
(i)\, Simplicial localizations of categories. {\em J. Pure and Appl.
Algebra} {\bf 17} (1980), 267-284.
\newline
(ii)\, Calculating simplicial localizations. {\em J. Pure and
Appl. Algebra} {\bf 18} (1980), 17-35.
\newline
(iii)\, Function complexes in homotopical algebra. {\em Topology}
{\bf 19} (1980), 427-440.
\newline
(iv)\, Equivalences between homotopy theories of diagrams.
{\em Algebraic Topology and Algebraic $K$-theory}, {\em Annals of Math. Studies}
{\bf 113}, Princeton University Press (1987), 180-205.
\bibitem{dhk}
W. Dwyer, P. Hirschhorn, D. Kan. Model categories and more general abstract
homotopy theory: a work in what we like to think of as progress. Preprint.


\bibitem{gz}
P. Gabriel, M. Zisman. {\em Calculus of fractions and homotopy theory.}
Ergebnisse der Math. {\bf 35}, Springer-Verlag, New York (1967).

\bibitem{Hirschhorn}
P. Hirschhorn. {\em Localization of model categories.} Book-preprint, available
at \newline {\tt http://www-math.mit.edu}.

\bibitem{descente}
A. Hirschowitz, C. Simpson. Descente pour les $n$-champs. Preprint
math/9807049


\bibitem{q}
D. Quillen.
{\em Homotopical algebra}. Springer Lecture Notes in Mathematics
{\bf 43} (1967).

\bibitem{Reedy}
C. Reedy. Homotopy theory of model categories. Preprint (1973) available from P.
Hirschhorn.


\bibitem{cmstruct}
C. Simpson.
A closed model structure for $n$-categories, internal $Hom$,
$n$-stacks and generalized Seifert-Van Kampen, alg-geom/9704006.

\bibitem{Tamsamani}
Z. Tamsamani. 
Sur des notions de $n$-cat\'egorie et $n$-groupoide non strictes 
via des ensembles multi-simpliciaux. 
{\em $K$-Theory} {\bf 16} (1999),  51-99; cf  alg-geom/9512006 and 9607010.
















\end{thebibliography}
\end{document}